\documentclass[12pt]{amsart}
\usepackage{geometry}                % 
\usepackage{amssymb}
\usepackage{amsmath}
\usepackage{amsthm}
\DeclareMathOperator{\cl}{cl}
\DeclareMathOperator{\Fix}{Fix}

\newtheorem{theorem}{Theorem}
\newtheorem*{lemma}{Lemma}
\newtheorem*{corollary}{Corollary}

\newtheorem*{remark}{Remark}

\title{Prevalence of Odometers in Cellular Automata}

\author{Ethan M. Coven}
\address{Department of Mathematics, Wesleyan University, Middletown CT 06457-0128}
\email{ecoven@wesleyan.edu}

\author{Marcus Pivato}
\address{Department of Mathematics,
Trent University,
Peterborough ON,
Canada K9L 1Z8}
\email{pivato@xaravve.trentu.ca}

\author{Reem Yassawi}
\address{Department of Mathematics,
Trent University,
Peterborough ON, Canada K9L 1Z8}
\email{ryassawi@trentu.ca}

\thanks{
This work was done in Spring~2005 while the second and third authors were van Vleck Visiting Professors of Mathematics at Wesleyan University.  The first author wishes to thank the lovely summer weather on Cape ~Cod for delaying the submission of this article.}

\subjclass[2000]{Primary 37B10, 37B15}

\keywords{odometer, cellular automaton}

\date{October 2005}

\begin{document}

\begin{abstract}
We consider left permutive cellular automata ~$\Phi$ with no memory and positive anticipation,  defined on the space of all doubly infinite sequences with entries from a finite alphabet.  For each such automaton  that is not one-to-one, there is a dense set of points ~$x$, which is large in another sense too, such that 
$\Phi : \cl \{\Phi^n(x) : n \ge 0\}
\to \cl \{\Phi^n(x) : n \ge 0\}$
is topologically conjugate to an odometer, the ``$+1$'' map on the countable product of finite cyclic groups.  We identify the odometer in several cases.
\end{abstract}

\maketitle

\section*{Introduction}

In this paper we show that for a certain class of one-dimensional cellular automata ~$\Phi$ defined on the space of all doubly infinite sequences with entries from a finite alphabet, there are many sequences ~$x$ such that
$$\Phi : \cl \{\Phi^n(x) : n \ge 0\}
\to \cl \{\Phi^n(x) : n \ge 0\}$$
is topologically conjugate to an odometer.  We then investigate the size of the set of such points.

We  use the most concrete definition of odometer.  
Let $(s_1,s_2,\dots)$ be a sequence of integers greater than $1$.  The \emph {$(s_1,s_2,\dots)$-adic odometer} is the ``$+1$'' map ~$\tau$ defined on the compact, abelian group
$$\mathbb{Z}(S) = \prod_{n\ge1}
\mathbb{Z}/s_n\mathbb{Z}$$
where $S = (s_1,s_2,\dots)$,
addition is ``with carrying''
and $\tau :\mathbb{Z}(S) \to \mathbb{Z}(S)$ is defined by
$$\tau(z) = z + (1,0,0,\dots).$$
 When $S$ is the constant sequence 
$(s,s,\dots)$,
$\tau :\mathbb{Z}(S) \to \mathbb{Z}(S)$
is called the \emph{$s$-adic odometer\/}
and denoted
$\tau :\mathbb{Z}(s) \to \mathbb{Z}(s)$.

Odometers are also called \emph{adding machines}.  They are characterized by being minimal, having uniformly equicontinuous powers, and having rational point spectrum with respect to Haar measure \cite{D}.

A \emph{cellular automaton} is a continuous, shift-commuting self-map, defined on the space  
of all doubly infinite sequences with entries from a finite alphabet.
  It is well-known that every cellular automaton $\Phi$ is given by a local rule~$\varphi$: for some $r \ge 0$, for all~$x$, and for all ~$i$, $-\infty < i < \infty$,
$$[\Phi(x)]_i = \varphi(x_{i-r},x_{i-r+1},\dots,x_{i+r}).$$
$\Phi$ has \emph{anticipation\/}
$r>0$ if and only if there exist
$\overline{t_{-r}},\overline{t_{-r+1}},\dots,\overline{t_{r-1}}$ such that
$\varphi(\overline{t_{-r}},\overline{t_{-r+1}},\dots,\overline{t_{r-1}},\cdot)$
is not the constant function.
$\Phi$ has \emph{no memory\/}
if and only if
$\varphi(t_{-r},t_{-r+1},\dots,t_r)$
depends only on $t_0,t_1,\dots,t_r$.
In this case we omit
$t_{-r},t_{-r+1},\dots,t_{-1}$
and write $\varphi(t_0,t_1,\dots,t_r)$. 
Finally $\Phi$, with no memory and anticipation~$r>0$, is
\emph{left permutive\/}
if and only if
for every $t_1,\dots,t_r$,
$\varphi(\cdot,t_1,t_2,\dots,t_r)$
is a permutation of the alphabet.

Composition of local rules is defined so that if $\varphi$ is the local rule of ~$\Phi$, then $\varphi^2$ is the local rule of ~$\Phi^2$, etc.  Thus if
$\varphi = \varphi(t_0,t_1,\dots,t_r)$, then
$$\varphi^2(t_0,t_1,\dots,t_{2r}) =
\varphi(\varphi(t_0,t_1,\dots,t_r),\varphi(t_1,t_2,\dots,t_{r+1}),\dots,
\varphi(t_r,\dots,t_{2r})).$$

Left permutive cellular automata are a well-studied class.  Among their pleasant properties are that they are finite-to-one,  map the space of all doubly infinite sequences onto itself, and   preserve Bernoulli $(\frac{1}{s},\frac{1}{s},\dots,\frac{1}{s})$-measure, where $s$ is the size of the alphabet.  See, for example, \cite{H,LM}.

The statements of the theorems and their proofs involve the 
\emph{one-sided cellular automaton\/}
~$\Phi_R$, defined on the space of all one-sided sequences by the same local rule as~$\Phi$.  

We show  for   left permutive cellular automata with no memory  that are not one-to-one, a trivial necessary condition  for
$\Phi : \cl \{\Phi^n(x) : n \ge 0\}
\to \cl \{\Phi^n(x) : n \ge 0\}$
to be topologically conjugate to an odometer is also sufficent.  We  show that the set of such points ~$x$ is  dense, and is large in another sense we make precise.  We  identify the odometer in a number of cases.

Note that odometers are one-to-one maps, while in most cases the cellular automata we consider are not.

\section{Existence of Odometers}

In this section we show 
for  left permutive cellular automata~$\Phi$ with no memory, 
a trivial necessary condition on ~$x$ with $(x_1,x_2,\dots)$
 $\Phi_R$-fixed for
$\Phi : \cl \{\Phi^n(x) : n \ge 0\}
\to \cl \{\Phi^n(x) : n \ge 0\}$
to be topologically conjugate to an odometer, namely that 
$\{\Phi^n(x): n \ge 0\}$ is infinite, is also sufficient.  

\begin{theorem}
Let $\Phi$ be a left permutive cellular automaton with no memory, defined of the space of all doubly infinite sequences with entries from a finite alphabet.  If

\begin{itemize} 
\item $\{\Phi^n(x) : n \ge 0\}$ is infinite and
\item $(x_1,x_2,\dots)$ is $\Phi_R$-fixed,  
\end{itemize} 

\noindent then
$\Phi : \cl \{\Phi^n(x) : n \ge 0\}
\to \cl \{\Phi^n(x) : n \ge 0\}$
is topologically conjugate to an odometer.
\end{theorem}

\begin{proof}

Let $k_1 $ be least such that
$(x_{-{k_1}},x_{-{k_1}+1},\dots)$
is not $\Phi_R$-fixed
and let $s_1$ be least such that
$(x_{-{k_1}},x_{-{k_1}+1},\dots)$
is  $\Phi_R^{s_1}$-fixed.
($s_1$ exists because $\Phi$ is left permutive.)
As a model of the inductive step,
let $k_2 > k_1$ be least such that
$(x_{-{k_2}},x_{-{k_2}+1},\dots)$
is not $\Phi_R^{s_1}$-fixed
and let  $s_2$ be least such that
$(x_{-{k_2}},x_{-{k_2}+1},\dots)$
is  $\Phi_R^{s_1s_2}$-fixed.
Continue.

We show that
$\Phi : \cl \{\Phi^n(x) : n \ge 0\}
\to \cl \{\Phi^n(x) : n \ge 0\}$
is topologically conjugate to the 
$(s_1,s_2,\dots)$-adic odometer.
To do this we show that the 
map of 
$\{\Phi^n(x) : n \ge 0\}$ 
into~$\mathbb{Z}(S)$,
$S = (s_1,s_2,\dots)$,
 defined by
$$\Phi^n(x) \mapsto \hbox{ ``base-$S$''
expansion of}~n,$$
and its inverse are uniformly continuous.
(The map is well-defined because
$\{\Phi^n(x) : n \ge 0\}$
is infinite.) 

For $i \ge 0$, the ``base-$S$''
expansions of~$m$ and of~$n$
agree at places $0,1,\dots,i$
(a measure of closeness in
~$\mathbb{Z}(S)$) if and only if
$\Phi^m(x)$ and $\Phi^n(x)$
agree at places $-k_i,-k_i+1,\dots,0$.
Since $(x_1,x_2,\dots)$  is  $\Phi$-fixed,
$\Phi^m(x)$ and $\Phi^n(x)$
agree at places $-k_i,-k_i+1,\dots,0$
if and only if they agree at places
$-k_i,-k_i+1,\dots,k_i$
(a measure of closeness in the space of all doubly infinite sequences).
Therefore the map is uniformly continuous on $\{\Phi^n(x) : n \ge 0\}$
and its inverse is uniformly continuous  on the image of this set.

Let $\Psi$ be the extension of this map to a homeomorphism defined on
$\cl \{\Phi^n(x) : n \ge 0\}$.
Since the image of 
$\{\Phi^n(x) : n \ge 0\}$,
the non-negative integers,
is dense in~$\mathbb{Z}(S)$,
$\Psi$ maps 
$\cl \{\Phi^n(x) : n \ge 0\}$
onto ~$\mathbb{Z}(S)$.
Since 
$\Psi \circ \Phi = \tau \circ \Psi$
on $\{\Phi^n(x) : n \ge 0\}$,
$\Psi \circ \Phi = \tau \circ \Psi$
on $\cl \{\Phi^n(x) : n \ge 0\}$
as well.  
\end{proof}

\begin{remark}
It follows from the proof that if  the odometer is the $(s_1,s_2,\dots)$-adic odometer and the alphabet has ~$s$ letters, then every $s_i \le s$.
\end{remark}

\begin{corollary} Let $\Phi$ be a left permutive cellular automaton with no memory, defined of the space of all doubly infinite sequences with entries from a finite alphabet.  If

\begin{itemize}
\item $\{\Phi^n(x) : n \ge 0\}$ is infinite, 
\item $(x_1,x_2,\dots)$ is $\Phi_R$-periodic with least period ~$q>1$, and 
\item $\Phi^q : \cl \{\Phi^{nq}(x) : n \ge 0\}
\to \cl \{\Phi^{nq}(x) : n \ge 0\}$
is topologically conjugate to the $(s_1,s_2,\dots)$-adic odometer, 
\end{itemize} 

\noindent then
$\Phi : \cl \{\Phi^n(x) : n \ge 0\}
\to \cl \{\Phi^n(x) : n \ge 0\}$
is topologically conjugate to the $(q,s_1,s_2,\dots)$-adic odometer. 
\end{corollary}

\begin{proof}
The topological conjugacy of
$\cl \{\Phi^{nq}(x) : n \ge 0\}$
to $\mathbb{Z}(s_1,s_2,\dots)$
is the extension of the map
$\Phi^{nq}(x) \mapsto \hbox{ ``base-$(s_1,s_2,\dots)$''
expansion of}~n.$
Extend this map to a map $\Psi_{0}$ of 
$\{\Phi^n(x) : n \ge 0\}$ into
$\mathbb{Z}(q,s_1,s_2,\dots)$
by $\Psi_0(\Phi^n(x)) :=  (v,w_1,w_2,\dots)$,
where $n = uq+v$, $0 \le v \le q-1$,
and $\Psi_{0}^{uq}(x) =(w_1,w_2,\dots)$.
Then $\Psi_0$ extends to a
topological conjugacy of
$\cl \{\Phi^n(x) : n \ge 0\}$ onto
$\mathbb{Z}(q,s_1,s_2,\dots)$.

\end{proof}

\section{Prevalence of Odometers}

In this section we identify  senses in which the set of points ~$x$ such that
$\Phi : \cl \{\Phi^n(x) : n \ge 0\}
\to \cl \{\Phi^n(x) : n \ge 0\}$
is  topologically conjugate to an 
odometer is large.

It is clear 
from looking at the form of the local rules of the powers of cellular automata
that every positive power of a left permutive cellular automaton with no memory is left permutive and has no memory.  This not the case with positive anticipation, although we do have the following.

\begin{lemma}\cite[Theorem 6.9]{H}
Let $\Phi$ be a left permutive
cellular automaton with no memory and positive anticipation, defined on the space of all doubly infinite sequences with entries from a 2-letter alphabet. Then $\Phi^n$ has positive anticipation for every $n \ge 1$.
\end{lemma}

The lemma is not true for larger alphabets, as shown by the following example.  Consider the  cellular automaton ~$\Phi$, defined on the space of all doubly infinite sequences with entries from
~$\{0,1,2\}$, and local rule
$\varphi(t_0, t_1) = t_0$,
except that 
$\varphi(0,1) = 2$ and
$\varphi(2,1) = 0$.
Then $\Phi^2$ is the identity map.

That this is essentially the only example is shown by the following

\begin{lemma}
Let $\Phi$ be a left permutive
cellular automaton with no memory, defined on the space of all doubly infinite sequences with entries from a finite alphabet.  
Then $\Phi^n$ has positive anticipation for every $n \ge 1$
if and only if $\Phi^m$ is not the identity map for every $m \ge 1$.
In particular, if $\Phi$ is not one-to-one, then $\Phi^n$ has positive anticipation for every $n \ge 1$.

\end{lemma}

\begin{proof}
Suppose that $\Phi^n$ has zero anticipation.  Since $\Phi^n$ is left permutive, it is a permutation of the alphabet.  Therefore $\Phi^{kn}$ is the identity map for some $k \ge 1$.
\end{proof}

\begin{theorem}
Let $\Phi$ be a  left permutive cellular automaton with no memory and positive anticipation, defined on the space of all doubly infinite sequences with entries from a finite alphabet.
Then exactly one of the following statements holds.

\begin{enumerate}

\item For every $\Phi_R$-fixed $(z_1,z_2,\dots)$, 
the set of points $x$ such that
$x_i = z_i$ for every  $i \ge 1$ and 
$\Phi : \cl \{\Phi^n(x) : n \ge 0\}
\to \cl \{\Phi^n(x) : n \ge 0\}$
is  topologically conjugate to an 
dometer is a dense $G_{\delta}$ subset of $\{x: x_i = z_i \mbox{ for every }i \ge 1  \}$.

\item For every $\Phi_R$-fixed $(z_1,z_2,\dots)$,
the set of points $x$
such that
$x_i = z_i$ for every  $i \ge 1$ and 
$\Phi : \cl \{\Phi^n(x) : n \ge 0\}
\to \cl \{\Phi^n(x) : n \ge 0\}$
is  topologically conjugate to an 
odometer is empty.

\end{enumerate}

If the alphabet has two letters or if the cellular automaton is not one-to-one, then (1) holds.
\end{theorem}

\begin{proof} 
It suffices to show that if $\Phi^n$ has positive anticipation for every $n \ge 1$, then (1) holds.

 Since $\{x: x_i = z_i \mbox{ for every }i \ge 1  \}$
is a complete metric space,
by Theorem~1 and its Corollary it is sufficient to show that the set of points in this set
with finite $\Phi$-orbits is a countable union of sets which are closed and nowhere dense in this set.

The set of points in
$\{x: x_i = z_i \mbox{ for every }i \ge 1  \}$ 
with finite $\Phi$-orbits is
$$\bigcup_{i\ge0}\bigcup_{j\ge0}
  (\Phi^{-i}(\Fix(\Phi^j) \cap
  \{x: x_i = z_i \mbox{ for every }i \ge 1  \}).$$         
For each $j \ge 1$,
$\Fix(\Phi^j) \cap \{x: x_i = z_i \mbox{ for every }i \ge 1  \}$
is a closed and nowhere dense subset of 
$\{x: x_i = z_i \mbox{ for every }i \ge 1  \}$.
It is nowhere dense because ~$\Phi^j$ has positive anticipation, and so any point in $\Fix(\Phi^j) \cap \{x: x_i = z_i \mbox{ for every }i \ge 1  \}$ can be changed 
arbitrarily far to the left so that it still is in $\{x: x_i = z_i \mbox{ for every }i \ge 1  \}$ but not in
$\Fix(\Phi^j)$.

Since every $\Phi^i$ is left permutive and hence a self-homeomorphism of 
$\{x: x_i = z_i \mbox{ for every }i \ge 1  \}$, 
each set in the double union above is a closed and nowhere dense subset of
$\{x: x_i = z_i \mbox{ for every }i \ge 1  \}$.

It follows from the two lemmas preceding the theorem that if the alphabet has two letters or if the cellular automaton is not one-to-one, then (1) holds.

\end{proof}

\begin{corollary}
Let $\Phi$ be a  left permutive cellular automaton with no memory and positive anticipation, defined on the space of all doubly infinite sequences with entries from a finite alphabet.
Then exactly one of the following statements holds.

\begin{enumerate}

\item For every $\Phi_R$-periodic  $(z_1,z_2,\dots)$,
 the set of points in 
$$\bigcup_{0\le k\le q-1} \{x: x_i = [\Phi^k(z)]_i \mbox{ for every }i \ge 1  \},$$ where  $q$ is the least period of $(z_1,z_2,\dots)$, such that
$\Phi : \cl \{\Phi^n(x) : n \ge 0\}
\to \cl \{\Phi^n(x) : n \ge 0\}$
is  topologically conjugate to an 
odometer is a dense $G_{\delta}$ subset of $$\bigcup_{0\le k\le q-1} \{x: x_i = [\Phi^k(z)]_i \mbox{ for every }i \ge 1  \}.$$

\item For every $\Phi_R$-periodic  $(z_1,z_2,\dots)$,
 the set of points in 
$$\bigcup_{0\le k\le q-1} \{x: x_i = [\Phi^k(z)]_i \mbox{ for every }i \ge 1  \},$$ where  $q$ is the least period of $(z_1,z_2,\dots)$, such that
$\Phi : \cl \{\Phi^n(x) : n \ge 0\}
\to \cl \{\Phi^n(x) : n \ge 0\}$
 is empty.
 
 \end{enumerate}

If the alphabet has two letters or if the cellular automaton is not one-to-one, then (1) holds.

\end{corollary}

\begin{proof}

Each set $\{x : x_i = [\Phi^k(z)]_i
\mbox{ for every } i \ge 1\}$, $k = 0,1,\dots,q-1$, contains a dense $G_\delta$ subset.  Since these sets are closed and pairwise disjoint, the union of the dense $G_\delta$ sets is a dense $G_\delta$ subset of $\bigcup_{0\le k\le q-1} \{x: x_i = [\Phi^k(z)]_i \mbox{ for every }i \ge 1  \}$.

\end{proof}

Another sense in which the set of points ~$x$ such that
$\Phi : \cl \{\Phi^n(x) : n \ge 0\}
\to \cl \{\Phi^n(x) : n \ge 0\}$
is  topologically conjugate to an 
odometer is large is given by the following.

\begin{theorem}
Let $\Phi$ be a  left permutive cellular automaton with no memory and positive anticipation, defined of the space of all doubly infinite sequences with entries from a finite alphabet.
Then exactly one of the following statements holds.

\begin{enumerate}

\item The set of points 
$x$ such that
$\Phi : \cl \{\Phi^n(x) : n \ge 0\}
\to \cl \{\Phi^n(x) : n \ge 0\}$
is  topologically conjugate to an 
odometer is  dense in the space of all doubly infinite sequences.

\item The set of points 
$x$ such that
$\Phi : \cl \{\Phi^n(x) : n \ge 0\}
\to \cl \{\Phi^n(x) : n \ge 0\}$
is  topologically conjugate to an 
odometer is  empty.

\end{enumerate}

If the alphabet has two letters or if the cellular automaton is not one-to-one, then (1) holds.

\end{theorem}

\begin{proof}
By \cite{BK} the set of $\Phi$-periodic points is dense in the space of all doubly infinite sequences.  
The result then follows from the elementary fact that if for every ~$i$,
$Y_i$ is dense in ~$X_i$, and if $\bigcup X_i$ is dense in~ $X$,
then $\bigcup Y_i$ is dense in~ $X$.

\end{proof}

\section{Identifying the Odometer}
In this section we identify the odometers in Theorem~1 and its Corollary for certain cellular automata. We assume, without loss of generality, that when  
~ $s$ is the size of the alphabet,  
the alphabet is ~$\mathbb{Z}/(s)$, the ring of integers modulo ~$s$, and we restrict our attention to left permutive cellular automata with no memory  and anticipation $r > 0$, whose local rules can be written in the form 
$t_0 + \theta(t_1,\dots,t_r)$.
Recall that
when the alphabet is ~$\mathbb{Z}/(2)$, every positive power of  such a cellular automaton with positive anticipation has positive anticipation.
It is easy to see that when the alphabet is ~$\mathbb{Z}/(2)$, the local rule of every left permutive cellular automaton with no memory must be of this form.  That this is not true for larger alphabets is shown by the example at the beginning of Section ~2.

\begin{theorem}
Let $\Phi$ be a left permutive cellular automaton with no memory and anticipation $r>0$, defined on the space of all doubly infinite sequences with entries from ~$\mathbb{Z}/(p)$, $p$ prime, and whose local rule is
$t_0 + \theta(t_1,\dots,t_r)$ for some function ~$\theta$.
If 

\begin{itemize}
\item $\{\Phi^n(x) : n \ge 0\}$ is infinite and 
\item $(x_1,x_2,\dots)$ is $\Phi_R$-fixed,
\end{itemize}  

\noindent then $\Phi : \cl \{\Phi^n(x) : n \ge 0\}
\to \cl \{\Phi^n(x) : n \ge 0\}$
is topologically conjugate to the $p$-adic odometer.
\end{theorem}

\begin{proof}
As in the proof of Theorem~1,
let $k_1 $ be least such that
$(x_{-{k_1}},x_{-{k_1}+1},\dots)$
is not  $\Phi_R$-fixed
and let $s_1$ be least such that
$(x_{-{k_1}},x_{-{k_1}+1},\dots)$
is   ~$\Phi_R^{s_1}$-fixed.
Since
$[\Phi_R^i(x_{-k_1},x_{-k_1+1},\dots)]_{-k_1} =
x_{-k_1} + i\theta(t_{-k_1+1},\dots,t_{-k_1+r})$
and  $\theta(t_{-k_1+1},\dots,t_{-k_1+r}) \ne 0$, $s_1 = p$.  Continuing as in the proof of Theorem~1, we find that
$s_2 = s_3 = \dots = p$.
\end{proof}

\begin{corollary}
Let $\Phi$ be a left permutive cellular automaton with no memory and anticipation $r>0$, defined on the space of all doubly infinite sequences with entries from ~$\mathbb{Z}/(p)$, $p$ prime, and whose local rule is
$t_0 + \theta(t_1,\dots,t_r)$
for some function ~$\theta$.
If 

\begin{itemize}
\item $\{\Phi^n(x) : n \ge 0\}$ is infinite and 
\item $(x_1,x_2,\dots)$ is $\Phi_R$-periodic with least period~$q>1$,
\end{itemize}  

\noindent then
$\Phi : \cl \{\Phi^n(x) : n \ge 0\}
\to \cl \{\Phi^n(x) : n \ge 0\}$
is topologically conjugate to the $(q,p,p,\dots)$-adic odometer.
\end{corollary}

\begin{theorem}
Let $\Phi$ be a left permutive cellular automaton with no memory and anticipation $r>0$, defined on the space of all doubly infinite sequences with entries from ~$\mathbb{Z}/(p^m)$, $p$ prime,  and whose local rule is
$t_0 + \theta(t_1,\dots,t_r)$ for some function ~$\theta$.
If 

\begin{itemize}
\item $\{\Phi^n(x) : n \ge 0\}$ is infinite and 
\item $(x_1,x_2,\dots)$ is $\Phi_R$-fixed,
\end{itemize}  

\noindent then
$\Phi : \cl \{\Phi^n(x) : n \ge 0\}
\to \cl \{\Phi^n(x) : n \ge 0\}$
is topologically conjugate to the $p$-adic odometer.
\end{theorem}

\begin{proof}
In the proof of Theorem~4,
the least $i > 0$ such that
$(x_{-{k_1}},x_{-{k_1}+1},\dots)$
is   $\Phi_R^i$-fixed is ~$p^{m_1}$ for some positive $m_1 \le m$.  Continuing this way, we get that
$\Phi : \cl \{\Phi^n(x) : n \ge 0\}
\to \cl \{\Phi^n(x) : n \ge 0\}$
is topologically conjugate to the 
$(p^{m_1},p^{m_2},\dots)$-adic odometer, which is topologically conjugate to the $p$-adic odometer. 
\end{proof}

\begin{corollary}
Let $\Phi$ be a left permutive cellular automaton with no memory and anticipation $r>0$, defined on the space of all doubly infinite sequences with entries from ~$\mathbb{Z}/(p^m)$, $p$ prime, and whose local rule is
$t_0 + \theta(t_1,\dots,t_r)$ for some function ~$\theta$.
If 

\begin{itemize}
\item $\{\Phi^n(x) : n \ge 0\}$ is infinite and
 
\item  $(x_1,x_2,\dots)$ is $\Phi_R$-periodic with least period~$q>1$,
\end{itemize}  

\noindent then
$\Phi : \cl \{\Phi^n(x) : n \ge 0\}
\to \cl \{\Phi^n(x) : n \ge 0\}$
is topologically conjugate to the $(q,p,p,\dots)$-adic odometer.
\end{corollary}

\begin{theorem}
Let $\Phi$ be a left permutive cellular automaton with no memory and anticipation $r>0$, defined on the space of all doubly infinite sequences with entries from ~$\mathbb{Z}/(p)$, $p$ prime, and whose local rule is
$at_0 + \theta(t_1,\dots,t_r)$ for some function ~$\theta$
and some $a \ne 0,1$.
If

\begin{itemize} 
\item $q$ is least such that $a^q = 1$,

\item $\{\Phi^n(x) : n \ge 0\}$ is infinite, and 

\item $(x_1,x_2,\dots)$ is $\Phi_R$-fixed,
\end{itemize}  

\noindent then
$\Phi : \cl \{\Phi^n(x) : n \ge 0\}
\to \cl \{\Phi^n(x) : n \ge 0\}$
is topologically conjugate to the 
$(q,p,p,\dots)$-adic odometer.
\end{theorem}

\begin{proof}
Use the proof of the Corollary to Theorem~1.
\end{proof} 

The reader is invited to state and prove for Theorem ~6 the result corresponding to the Corollary to Theorem~4.

\end{document}